\documentclass{tac}

\usepackage{amsmath}
\usepackage{amssymb}
\usepackage{enumitem}

\usepackage[all,cmtip]{xy}

\input diagxy

\usepackage[colorlinks=true]{hyperref}
\hypersetup{allcolors=[rgb]{0.1,0.1,0.4}}

\author{Xiaoye Tang}

\thanks{The author acknowledges the support of National Natural Science Foundation of China (No. 12071319).}

\address{School of Mathematics, Sichuan University\\
 Chengdu 610064, China\\
}

\title{Completeness of quantaloid-enriched categories up to Morita equivalence}

\copyrightyear{2025}


\eaddress{tang.xiaoye@qq.com}

\newtheorem{thm}{Theorem}

\newtheorem{lem}{Lemma}
\newtheorem{prop}{Proposition}

\newtheoremrm{rem}{Remark}
\newtheoremrm{defn}{Definition}
\newtheoremrm{exmp}{Example}
\newtheoremrm{ques}{Question}

\mathrmdef{Hom}

\DeclareMathOperator{\ub}{ub}
\DeclareMathOperator{\lb}{lb}

\def\oto{{\bfig\morphism<180,0>[\mkern-4mu`\mkern-4mu;]\place(86,0)[\circ]\efig}}

\def\nra{\relbar\joinrel\joinrel\mapstochar\joinrel\joinrel\rightarrow}
\def\nla{\leftarrow\joinrel\joinrel\joinrel\mapstochar\joinrel\relbar}

\newcommand{\ra}{\rightarrow}
\newcommand{\la}{\leftarrow}

\newcommand{\lda}{\swarrow}
\newcommand{\rda}{\searrow}

\newcommand{\rat}{\rightarrowtail}
\newcommand{\bv}{\bigvee}
\newcommand{\bw}{\bigwedge}

\newcommand{\dv}{\dashv}

\renewcommand{\phi}{\varphi}

\newcommand{\lam}{\lambda}

\newcommand{\CC}{\mathcal{C}}

\newcommand{\CM}{\mathcal{M}}

\newcommand{\CP}{\mathcal{P}}
\newcommand{\CQ}{\mathcal{Q}}

\newcommand{\sY}{\mathsf{Y}}

\newcommand{\sm}{\mathsf{m}}

\newcommand{\sfs}{\mathsf{s}}

\newcommand{\bbA}{\mathbb{A}}
\newcommand{\bbB}{\mathbb{B}}
\newcommand{\bbC}{\mathbb{C}}

\newcommand{\Cat}{{\bf Cat}}

\newcommand{\Dist}{{\bf Dist}}
\newcommand{\Inf}{{\bf Inf}}
\newcommand{\Map}
{\mathbf{Map}}

\newcommand{\Sup}{{\bf Sup}}

\newcommand{\QCat}{\CQ\text{-}\Cat}

\newcommand{\QMC}{\CQ\text{-}\bf MC}
\newcommand{\QMCC}{\CQ\text{-}\bf MCC}

\newcommand{\QDist}{\CQ\text{-}\Dist}

\newcommand{\QInf}{\CQ\text{-}\Inf}

\newcommand{\QSup}{\CQ\text{-}\Sup}

\newcommand{\CPd}{\CP^{\dag}}

\newcommand{\sYd}{\sY^{\dag}}
\newcommand{\co}{{\rm co}}
\newcommand{\op}{{\rm op}}

\renewcommand{\leq}{\leqslant}
\renewcommand{\geq}{\geqslant}

\numberwithin{equation}{section}

\begin{document}

\maketitle
\begin{abstract}
For a small quantaloid $\mathcal{Q}$, we introduce $\mathcal{M}$-(co)complete $\mathcal{Q}$-categories, i.e., (co)complete $\mathcal{Q}$-categories up to Morita equivalence, as Eilenberg--Moore algebras of the presheaf monad on the category of $\mathcal{Q}$-categories and left adjoint $\mathcal{Q}$-distributors, and characterize such $\mathcal{Q}$-categories through $\mathcal{M}$-(co)tensoredness and $\mathcal{M}$-conical (co)completeness.
\end{abstract}


\section{Introduction}
 
Morita equivalence is an important concept in category theory, originally introduced by Morita \cite{Morita1958} (also see \cite{Bass1968}) in the context of the category of modules. Since then, the notion has been generalized to more abstract settings \cite{FisherPalmquist1975,Newell1972,Polin1974}, including enriched categories, and can be expressed in terms of profunctors or distributors, as observed  by Lawvere and B\'{e}nabou \cite{Lawvere1973,Benabou1973}.

This paper concerns Morita equivalences between categories enriched over a quantaloid. A quantaloid $\CQ$ \cite{Rosenthal1996,Stubbe2005,Stubbe2014} is a category enriched in the category $\Sup$ of complete lattices and sup-preserving functions. Two $\CQ$-categories $\bbA$ and $\bbB$ are said to be Morita equivalent \cite{Kelly1982,Stubbe2005} if they satisfy any of the following equivalent conditions:
\begin{itemize}
\item $\bbA\cong\bbB$ in $\QDist$;
\item $\CP\bbA\cong\CP\bbB$ in $\QCat$,
\item $\CPd\bbA\cong\CPd\bbB$ in $\QCat$,
\item $\bbA_{cc}\cong\bbB_{cc}$ in $\QCat$, 
\end{itemize}
where $\CP\bbA$, $\CPd\bbA$ and $\bbA_{cc}$ denote the presheaf category, the copresheaf category and the Cauchy completion of $\bbA$, respectively.

Therefore, a $\CQ$-category can be viewed as a presentation of its Cauchy completion, and two $\CQ$-categories are Morita equivalent if they have the same Cauchy completion. When one considers $\CQ$-categories up to Morita equivalence, one may equivalently replace $\QCat$ by the category $\Map(\QDist)$ of $\CQ$-categories and left adjoint $\CQ$-distributors, because left adjoint $\CQ$-distributors between Cauchy complete $\CQ$-categories are essentially $\CQ$-functors.  

It is well known that the Eilenberg--Moore algebras for the presheaf monad
\[(\CP,\sfs,\sY)\]
on the $2$-category $\QCat$ are precisely the skeletal cocomplete $\CQ$-categories \cite{Heymans2010}. Motivated by this correspondence, we seek to formulate notions of completeness and cocompleteness up to Morita equivalence by establishing the presheaf and copresheaf monads on $\Map(\QDist)$. In Section \ref{presheaf and copresheaf monads}, we introduce the monads 
\[(\widehat{\CP},\sm,\iota)\quad\text{and}\quad (\widehat{\CPd},\sm^\dagger,\iota^\dagger),\]
and prove the following results (Proposition \ref{Pre-md-al} and \ref{coPre-md-al}), which establish the structures of their Eilenberg--Moore categories:
\begin{itemize}
\item Let $\bbA$ and $\bbB$ be $\widehat{\CP}$-algebras. A left adjoint $\CQ$-distributor $\zeta\colon \bbA\oto\bbB$ is the $\widehat{\CP}$-algebra homomorphism if and only if it is a right adjoint in $\QDist$.
\item Let $\bbA$ and $\bbB$ be $\widehat{\CPd}$-algebras. A left adjoint $\CQ$-distributor $\zeta\colon \bbA\oto\bbB$ is the $\widehat{\CPd}$-algebra homomorphism if and only if $\zeta^*$ is a left adjoint in $\QDist$.
\end{itemize}
In Section \ref{M-completeness}, we give formal definitions of cocompleteness and completeness up to Morita equivalence and further show that, as in the classical theory, cocompleteness up to Morita equivalence is characterized by being tensored and conically cocomplete up to Morita equivalence. 

\section{Quantaloid-enriched categories}\label{Quantaloid}

This section reviews essential concepts and results on quantaloid-enriched categories that will be required in the following discussion. 

A \emph{quantaloid} \cite{Rosenthal1996,Stubbe2005,Stubbe2014} $\CQ$ is a category enriched in the symmetric monoidal closed category $\Sup$. Explicitly, a quantaloid $\CQ$ is a (possibly large) category with a class $\CQ_0$ of objects, such that
\begin{itemize}
\item each hom-set $\CQ(X,Y)$ ($X,Y\in\CQ_0$) is a complete lattice, and
\item the composition $\circ$ of $\CQ$-arrows preserves joins in each variable, i.e.
\[g\circ\Big(\bv\limits_{i}f_i\Big)=\bv\limits_{i}(g\circ f_i)\quad\text{and}\quad \Big(\bv\limits_{j} g_j\Big)\circ f=\bv\limits_{j}(g_j\circ f)\]
for all $\CQ$-arrows $f,f_i\colon X\to Y$ and $g,g_i\colon Y\to Z$ ($i\in I$).
\end{itemize}
In what follows, $\CQ$ always denotes a \emph{small}    quantaloid with a set $\CQ_0$ and a set $\CQ_1$ of $\CQ$-arrows. We denote the top and the bottom element of $\CQ(X,Y)$ by $\top_{X,Y}$ and $\bot_{X,Y}$,  respectively, and the identity arrow on $X\in\CQ_0$ by $1_X$.

Given $\CQ$-arrows $f\colon X\to Y$, $g\colon Y\to Z$, $h\colon X\to Z$, the corresponding adjoints induced by the compositions
\[-\circ f\dv -\swarrow f\colon\CQ(X,Z)\to\CQ(Y,Z),\]
\[g\circ-\dv g\searrow -\colon\CQ(X,Z)\to\CQ(X,Y)\]
satisfy
\[g\circ f\leq h\iff g\leq h\swarrow f\iff f\leq g\searrow h,\]
where the operations $\swarrow$, $\searrow$ are called \emph{left} and \emph{right implications} in $\CQ$, respectively. 

An \emph{adjunction} in a quantaloid $\CQ$ is a pair of $\CQ$-arrows $f\colon X\to Y$ and $g\colon Y\to X$ in $\CQ$, denoted by $f\dv g\colon Y\to X$, such that
\[1_X\leq g\circ f\quad\text{and}\quad f\circ g\leq 1_Y.\]
In this case, $f$ is a left adjoint of $g$ and $g$ is a right adjoint of $f$. 

\begin{prop}\label{Q-dv-cd} \cite{Heymans2010}
Let $f\dv g$ in a quantaloid $\CQ$. Then, for all $h,h'\in\CQ_1$, such that the corresponding operations are defined, the following identities hold:
\begin{enumerate}[label=(\arabic*)]
\item \label{Q-dv-cd-1} $h\circ f=h\lda g$, $g\circ h=f\rda h$.
\item \label{Q-dv-cd-3} 
$(h\rda h')\circ f=h\rda(h'\circ f)$, $g\circ(h'\lda h)=(g\circ h')\lda h$.
\end{enumerate}
\end{prop}

\begin{prop} \label{ad-Q}
If $f\dv g\colon Y\to X$ in a quantaloid $\CQ$, then
\begin{equation}
g=f\rda 1_Y\quad\text{and}\quad f=1_Y\lda g.
\end{equation}
\end{prop}

Therefore, the right adjoint of a morphism in $\CQ$, when it exists, is necessarily unique. In what follows we define 
\begin{equation} \label{u*-def}
f^*:=f\rda 1_Y\colon Y\to X
\end{equation}
for each morphism $f\colon X\to Y$ in $\CQ$. $f$ is called a \emph{map} in $\CQ$ \cite{Heymans2010} if $f\dv f^*$, and we denote by
\[\Map(\CQ)\]
the subcategory of $\CQ$ whose objects are the same as $\CQ$, and whose morphisms are maps in $\CQ$.

Since it always holds that $f\circ f^*=f\circ(f\rda 1_Y)\leq 1_Y$, we immediately obtain the following lemma:

\begin{lem} \label{u-map-leq}
A morphism $f\colon X\to Y$ is a map in $\CQ$ if, and only if, $1_X\leq f^*\circ f$.
\end{lem}

By Proposition~\ref{ad-Q} together with Lemma~\ref{u-map-leq}, a morphism $f\colon X\to Y$ in $\CQ$ is a right adjoint if and only if $1_Y\leq f\circ(1_X\lda f)$.

A \emph{$\CQ$-typed set} is a set $A$ equipped with a function $\left|-\right| \colon A\to\CQ_0$ assigning to each element $x\in A$ its type $|x|\in\CQ_0$. A \emph{$\CQ$-category} $\bbA$ consists of a $\CQ$-typed set $\bbA_0$ and hom-arrows $\bbA(x,y)\in\CQ(| x|,|y|)$ such that
\begin{enumerate}[label=(\arabic*)]
\item $1_{| x|}\leq \bbA(x,x)$ for all $x\in \bbA_0$;
\item $\bbA(y,z)\circ\bbA(x,y)\leq \bbA(x,z)$ for all $x,y,z\in\bbA_0$.
\end{enumerate}

Given a $\CQ$-category $\bbA$, there is a natural underlying preorder $\leq$ on $\bbA_0$ given by 
\[x\leq y\ \iff \  |x|=|y|=X\ \text{and}\ 1_X\leq\bbA(x,y).\]
For each $X\in\CQ_0$, the objects in $\bbA$ with type $X$ constitute a subset of the underlying preordered set $\bbA_0$ and we denote it by $\bbA_X$. Two objects $x,y$ in $\bbA$ are \emph{isomorphic} if $x\leq y$ and $y\leq x$, written $x\cong y$. $\bbA$ is \emph{skeletal} if no two different objects in $\bbA$ are isomorphic.

A \emph{$\CQ$-functor} $F\colon \bbA\to \bbB$ between $\CQ$-categories is a function $F\colon \bbA_0\to \bbB_0$ such that 
\begin{enumerate}[label=(\arabic*)]
\item $F$ is type-preserving, i.e. $|x|=|F(x)|$ for all $x\in\bbA_0$;
\item $\bbA(x,y)\leq \bbB(F(x),F(y))$ for all $x,y\in \bbA_0$.
\end{enumerate}
With the pointwise (pre)order of $\CQ$-functors given by
\[F\leq G\iff \forall x\in \bbA_0:\ F(x)\leq G(x),\]
$\CQ$-categories and $\CQ$-functors constitute a $2$-category
\[\QCat.\]
A pair of $\CQ$-functors $F\colon\bbA\to\bbB$ and $G\colon\bbB\to\bbA$ forms an \emph{adjunction} in $\QCat$, written $F\dashv G\colon \bbB\to\bbA$, if 
\[1_\bbA\leq G\circ F \quad\text{and}\quad F\circ G\leq 1_\bbB,\] 
where $1_\bbA$ and $1_\bbB$ denote the identity $\CQ$-functors on $\bbA$ and $\bbB$, respectively. In this case, $F$ is said to be a \emph{left adjoint} of $G$ and $G$ a \emph{right adjoint} of $F$.

A \emph{$\CQ$-distributor} $\phi\colon\bbA\oto\bbB$ is a 
function $\bbA_0\times \bbB_0\to\CQ_1$ assigning to each pair $(x,y)\in\bbA_0\times\bbB_0$ a $\CQ$-arrow $\phi(x,y)\in\CQ(|x|,|y|)$, such that
\begin{enumerate}[label=(\arabic*)]
\item $\forall x\in \bbA_0$, $\forall y,y'\in \bbB_0$, $\bbB(y',y)\circ\phi(x,y')\leq\phi(x,y)$;
\item $\forall x,x'\in \bbA_0$, $\forall y\in \bbB_0$, $\phi(x',y)\circ \bbA(x,x')\leq\phi(x,y)$.
\end{enumerate}
$\CQ$-categories and $\CQ$-distributors constitute a quantaloid
\[\QDist\]
in which 
\begin{itemize}
\item the local order is defined pointwise: for $\CQ$-distributors $\phi,\psi\colon \bbA\oto \bbB$,
\[\phi\leq\psi\ \iff \ \forall x\in \bbA_0,\ \forall y\in\bbB_0,\ \phi(x,y)\leq \psi(x,y);\]
\item the composition $\psi\circ\phi\colon \bbA\oto\bbC$ of $\CQ$-distributors $\phi\colon\bbA\oto\bbB$ and $\psi\colon \bbB\oto\bbC$ is given by
\[\forall x\in\bbA_0,\ \forall z\in\bbC_0,\ (\psi\circ\phi)(x,z)=\bv\limits_{y\in \bbB_0}\psi(y,z)\circ\phi(x,y);\]
\item the identity $\CQ$-distributor on a $\CQ$-category $\bbA$ is the hom-arrows of $\bbA$ and will be denoted by $\bbA\colon \bbA\oto \bbA$;
\item for $\CQ$-distributors $\phi\colon \bbA\oto \bbB$, $\psi\colon \bbB\oto \bbC$, $\eta\colon \bbA\oto \bbC$, the left implication $\eta\swarrow\phi\colon \bbB\oto \bbC$ and the right implication $\psi\searrow\eta\colon \bbA\oto\bbB$ are given by
\[\forall y\in\bbB_0,\ \forall z\in\bbC_0,\ (\eta\lda\phi)(y,z)=\bw\limits_{x\in \bbA_0}\eta(x,z)\lda\phi(x,y)\]
and
\[\forall x\in\bbA_0,\ \forall y\in\bbB_0,\ (\psi\rda\eta)(x,y)=\bw\limits_{z\in \bbC_0}\psi(y,z)\rda\eta(x,z).\]
\end{itemize}

An adjunction $\phi\dv\psi\colon \bbB\oto\bbA$ in the quantaloid $\QDist$ consists of a pair of $\CQ$-distributors $\phi\colon\bbA\oto\bbB$ and $\psi\colon\bbB\oto\bbA$ such that 
\[\bbA\leq\psi\circ\phi\quad\text{and}\quad \phi\circ\psi\leq\bbB.\]
Every $\CQ$-functor $F\colon \bbA\to \bbB$ induces an adjunction $F_\natural\dashv F^\natural\colon \bbA\rightharpoonup \bbB$ in $\QDist$, where $F_\natural$ and $F^\natural$, called the \emph{graph} and \emph{cograph} of $F$, are defined by   
\[F_\natural(x,y)=Y(F(x),y)\quad\text{and}\quad F^\natural(y,x)=Y(y,F(x))\]
for all $x\in \bbA_0$ and $y\in \bbB_0$. This correspondence allows adjunctions in $\QCat$ to be equivalently expressed in terms of the graphs and cographs of the functors, as summarized in the following Proposition.

\begin{prop} \label{ad-f-d} \cite{Stubbe2003}
Let $F\colon \bbA\to \bbB$ and $G\colon \bbB\to \bbA$ be a pair of $\CQ$-functors. The following conditions are equivalent:
\begin{enumerate}[label=(\arabic*)]
\item $F\dashv G\colon \bbB\to\bbA$.
\item\label{ad-f-d-rd}$F_\natural=G^\natural$.
\item\label{ad-f-d-cg}$ G^\natural\dashv F^\natural\colon \bbA\oto \bbB$.
\item \label{ad-f-d-g}$G_\natural\dv F_\natural\colon\bbB\oto\bbA$.
\end{enumerate}
\end{prop}

The assignment sending each $\CQ$-functor $F\colon \bbA \to \bbB$ to the $\CQ$-distributor $F_\natural\colon \bbA \oto \bbB$ determines a $2$-functor 
\[(-)_\natural\colon\QCat \to \Map(\QDist)\] which is identity on objects.

\begin{rem}
The dual of a $\CQ$-category $\bbA$ is a $\CQ^{\op}$-category $\bbA^{\op}$, defined by $\bbA^{\op}_0=\bbA_0$ and $\bbA^{\op}(x,y)=\bbA(y,x)$ for all $x,y\in\bbA_0$.
Each $\CQ$-functor $F\colon\bbA\to\bbB$ induces a $\CQ^{\op}$-functor $F^{\op}\colon\bbA^{\op}\to\bbB^{\op}$ with the same action on objects, satisfying $(F')^{\op}\leq F^{\op}$ whenever $F\leq F'$.
Similarly, each $\CQ$-distributor $\phi\colon\bbA\oto\bbB$ corresponds bijectively to a $\CQ^{\op}$-distributor $\phi^{\op}\colon\bbB^{\op}\oto\bbA^{\op}$ defined by $\phi^{\op}(y,x)=\phi(x,y)$ for all $x\in\bbA_0$, $y\in\bbB_0$.
Consequently, as already observed in \cite{Stubbe2005}, one obtains a 2-isomorphism
\begin{equation}\label{Q-Cat-op}
(-)^\op\colon\QCat\cong(\CQ^\op\text{-}\Cat)^\co
\end{equation}
together with an isomorphism of quantaloids 
\begin{equation}\label{Q-Dist-op}
(-)^\op\colon \QDist\cong(\CQ^\op\text{-}\Dist)^\op.
\end{equation}
\end{rem}

For each $X\in\CQ_0$, write $*_X$ for the $\CQ$-category with only one object $*$ of type $|*|=X$ and hom-arrow $1_X$.
A \emph{presheaf} on a $\CQ$-category $\bbA$ is a $\CQ$-distributor $\mu\colon \bbA \oto *_X$ with $X\in\CQ_0$. Presheaves on a $\CQ$-category $\bbA$ constitute a $\CQ$-category $\CP\bbA$ in which
\[|\mu|=X\quad\text{and}\quad \CP\bbA(\mu,\mu')=\mu'\lda \mu\]
for all $\mu\colon\bbA\oto *_X$ and $\mu'\colon\bbA\oto *_Y$.
Dually, a \emph{copresheaf} on a $\CQ$-category $\bbA$ is a $\CQ$-distributor $\lam\colon *_X\oto \bbA$. Copresheaves on $\bbA$ constitute a $\CQ$-category $\CPd \bbA$ in which
\[|\lam|=X\quad\text{and}\quad \CPd\bbA(\lam,\lam')=\lam'\rda \lam\]
for all $\lam\colon*_X\oto\bbA$ and $\lam'\colon *_Y\oto\bbA$. It is straightforward to verify that $\CPd\bbA\cong(\CP({\bbA^\op}))^\op$. In particular, we denote $\CP(*_X)=\CP X$ and $\CPd(*_X)=\CPd X$ for each $X\in\CQ_0$. 

For each $\CQ$-category $\bbA$, the \emph{Yoneda embedding} (\emph{co-Yoneda embedding}) is the $\CQ$-functor 
\[\sY_\bbA\colon\bbA\to\CP\bbA,\quad a\mapsto\bbA(-,a)\quad(\text{resp.}\ \sYd_\bbA\colon\bbA\to\CPd,\quad a\mapsto\bbA(a,-)).\]
The following statement is the well known \emph{Yoneda Lemma}.

\begin{lem}\label{Yoneda}(Yoneda). \cite{Stubbe2005}
For all $a\in\bbA_0$, $\mu\in(\CP\bbA)_0$ and $\lambda\in(\CPd\bbA)_0$,
\begin{equation}\label{Yoneda-eq}
\CP\bbA(\sY_\bbA (a),\mu)=\mu(a)\quad\text{and}\quad \CPd\bbA(\lambda,\sYd_\bbA (a))=\lambda(a).
\end{equation}
\end{lem}

A $\CQ$-category $\bbA$ is \emph{cocomplete} if the \emph{Yoneda embedding} $\sY_\bbA\colon\bbA\to\CP\bbA$ has a left adjoint $\sup\nolimits_\bbA\colon\CP\bbA\to\bbA$ in $\QCat$; that is,
\begin{equation}\label{sup}
\bbA(\sup\nolimits_\bbA(\mu),-)=\CP\bbA(\mu,\sY_\bbA(-))=\bbA\swarrow\mu
\end{equation}
for all $\mu\in(\CP\bbA)_0$. Dually, a $\CQ$-category $\bbA$ is \emph{complete} if the \emph{co-Yoneda embedding} $\sYd_\bbA$ has a right adjoint $\inf\nolimits_\bbA\colon\CPd \bbA\to\bbA$ in $\QCat$; that is,
\begin{equation}\label{inf}
\bbA(-,\inf\nolimits_\bbA (\lambda))=\CPd\bbA(\sYd_\bbA (-),\lambda)=\lambda\searrow\bbA
\end{equation}
for all $\lambda\in(\CPd \bbA)_0$. 

\begin{exmp} Let $\bbA$ be a $\CQ$-category.
\begin{enumerate}[label=(\arabic*)]
\item $\CP\bbA$ is a complete $\CQ$-category in which
\[\sup\nolimits_{\CP\bbA}(\Phi)=\Phi\circ(\sY_\bbA)_\natural\quad\text{and}\quad\inf\nolimits_{\CP\bbA}(\Psi)=\Psi\rda(\sY_\bbA)_\natural\]
for all $\Phi\in(\CP\CP\bbA)_0$ and $\Psi\in(\CPd\CP\bbA)_0$. 
\item $\CPd\bbA$ is a complete $\CQ$-category in which
\[\sup\nolimits_{\CPd\bbA}(\Phi)=(\sYd_\bbA)^\natural\lda\Phi\quad\text{and}\quad\inf\nolimits_{\CPd\bbA}(\Psi)=(\sYd_\bbA)^\natural\circ\Psi\]
for all $\Phi\in(\CP\CPd\bbA)_0$ and $\Psi\in(\CPd\CPd\bbA)_0$.
\end{enumerate}
\end{exmp}

Every $\CQ$-functor $F\colon \bbA\to\bbB$ induces two adjunctions\[F^\ra\dv F^\la\colon\CP\bbB\to\CP\bbA\quad\text{and}\quad F^{\nla}\dv F^{\nra}\colon\CPd\bbA\to\CPd\bbB\] defined by
\[F^\ra(\mu)=\mu\circ F^\natural,\quad F^\la(\lam)=\lam\circ F_\natural\]
and
\[ F^{\nla}(\lam')=F^\natural\circ\lam',\quad F^{\nra}(\mu')=F_\natural\circ\mu',\]
respectively.
The assignments
\[F\colon \bbA\to\bbB\ \mapsto \ \CP F=F^\ra\colon\CP\bbA\to\CP\bbB\]
and
\[F\colon\bbA\to\bbB\ \mapsto\ \CPd F=F^{\nra}\colon\CPd\bbA\to\CPd\bbB\]
yield the endofunctors $\CP$ and $\CPd$ on $\QCat$. The functor $\CP$ can be made into a monad $(\CP,\sfs,\sY)$, called the \emph{presheaf monad}, with unit given by the Yoneda embedding $\sY_\bbA\colon\bbA\to\CP\bbA$ and multiplication given by
\[\sfs_\bbA\colon\CP\CP\bbA\to\CP\bbA,\quad \sfs_\bbA(\Phi)=\sup\nolimits_{\CP\bbA}(\Phi)=\Phi\circ(\sY_\bbA)_\natural.\]
The presheaf monad $(\CP,\sfs,\sY)$ is a KZ-doctrine, or KZ-monad. The precise notion is as follows: A monad $(T,m,e)$ on a locally ordered category $\CC$ is a \emph{KZ-monad} (resp. \emph{co-KZ-monad}) \cite{Hofmann2014,Kock1995,Zoeberlein1976} if $T$ is a $2$-functor, and for all object $X$ in $\CC$, there is a string of adjunctions
\[Te_X\dv m_X\dv e_{TX}\quad(\text{resp.}\ e_{TX}\dv m_X\dv Te_X,).\]
A $T$-algebra for a KZ-monad (resp. co-KZ-monad) is a pair $(X,h)$  with  $h\circ e_X=1_X$, in which case $h\dv e_X$ (resp. $e_X\dv h$). Consequently, a $\CQ$-category $\bbA$ is a $\CP$-algebra if and only if $\sY_\bbA\colon\bbA\to\CP\bbA$ has a left inverse (hence $\bbA$ is skeletal), and in this case the left inverse is necessarily a left adjoint of $\sY_\bbA$. A $\CP$-algebra homomorphism $F\colon\bbA\to\bbB$ is a $\CQ$-functor $f\colon\bbA\to\bbB$ such that
\[\sup\nolimits_\bbB\circ\CP F=F\circ\sup\nolimits_\bbA,\]
which is equivalent to $F$ being a left adjoint $\CQ$-functor. Therefore, the category of $\CP$-algebras and $\CP$-algebra homomorphisms is just the category
\[\QSup\]
of skeletal cocomplete $\CQ$-categories and left adjoints $\CQ$-functors. Similarly, the functor $\CPd$ can also be made into a monad $(\CPd,\sfs^\dagger,\sYd)$, called the \emph{copresheaf monad}, with the unit $\sYd$ given by the co-Yoneda embedding and the multiplication $\sfs^\dag$ given by 
\[\sfs^\dag_\bbA\colon\CPd\CPd\bbA\to\CPd\bbA,\quad \sfs^\dag_\bbA(\Psi)=\inf\nolimits_{\CPd\bbA}(\Psi)=(\sYd_\bbA)^\natural\circ\Psi.\]
The copresheaf monad $(\CPd,\sfs^\dag,\sYd)$ is a co-KZ monad, thus the category of $\CPd$-algebras and $\CPd$-algebra homomorphisms is the category
\[\QInf\]
of skeletal complete $\CQ$-categories and right adjoints $\CQ$-functors. Indeed, every presheaf $\mu$ and every copresheaf $\lambda$ on $\bbA$ is, respectively, the supremum of $(\sY_\bbA)^\ra(\mu)$ and the infimum of $(\sYd_\bbA)^{\nra}(\lambda)$, that is,
\[\mu=\sup\nolimits_{\CP\bbA}({\sY_\bbA}^\ra(\mu)),\quad \lam=\inf\nolimits_{\CPd\bbA}({\sYd_\bbA}^{\nra}(\lam)).\]

In a $\CQ$-category $\bbA$, the \emph{tensor} (resp. \emph{cotensor}) of $a\in\bbA_0$ and $f\in(\CP |a|)_0$ (resp. $g\in(\CPd|a|)_0$), when it exists, is an object $f\otimes a$ (resp. $g\rat a$) in $\bbA_0$ of the type $|f\otimes a|=|f|$ (resp. $|g\rat a|=|g|$) such that
\begin{equation}\label{tensor-cotensor}
\bbA(f\otimes a,-)=\bbA(a,-)\swarrow f\quad(\text{resp.}\ \bbA(-,g\rat a)=g\searrow\bbA(-,a)).
\end{equation}
A $\CQ$-category $\bbA$ is \emph{tensored} (resp. \emph{cotensored}) if $f\otimes a$ (resp. $f\rat a$) exists for all $a\in\bbA_0$ and $f\in({\CP|a|})_0$ (resp. $f\in({\CP|a|})_0$). Indeed, in a complete $\CQ$-category $\bbA$, (\ref{tensor-cotensor}) indicates that there are adjunctions in $\QCat$
\[(-\otimes a)\dv\bbA(a,-)\colon\bbA\to\CP|a|\quad\text{and}\quad \bbA(-,a)\dv(-\rat a)\colon\CPd|a|\to \bbA\]
for all $a\in\bbA_0$.

A left adjoint $\CQ$-distributor $\phi\colon\bbA\oto\bbB$ is said to \emph{converge} if there exists a $\CQ$-functor $F\colon \bbA\to\bbB$ such that $\phi=F_\natural$.
A $\CQ$-category $\bbA$ is \emph{Cauchy complete}, if for any $\CQ$-category $\bbB$, every left adjoint $\CQ$-distributor $\phi\colon\bbB\oto\bbA$ converges.

For any $\CQ$-category $\bbA$, we define the \emph{Cauchy completion} $\bbA_{cc}$ of $\bbA$ by
\[(\bbA_{cc})_0=\{\mu\in(\CP\bbA)_0\ |\ \mu\ \text{is a left adjoint presheaf} \}\]
which is the full subcategory of $\CP\bbA$.

\begin{prop}\label{Morita-eq}\cite{Stubbe2005}
For $\CQ$-categories $\bbA$ and $\bbB$, the following are equivalent:
\begin{enumerate}[label=(\arabic*)]
\item $\bbA\cong\bbB$ in $\QDist$;
\item $\bbA_{cc}\cong\bbB_{cc}$ in $\QCat$;
\item $\CP\bbA\cong\CP\bbB$ in $\QCat$;
\item $\CPd\bbA\cong\CPd\bbB$ in $\QCat$.
\end{enumerate}
\end{prop}

We say that $\CQ$-categories $\bbA$ and $\bbB$ are \emph{Morita equivalent} \cite{Stubbe2005} if they satisfy the equivalent conditions in Proposition \ref{Morita-eq}. Conceptually, a $\CQ$-category can be understood  as a presentation of its Cauchy completion, and if one is concerned with $\CQ$-categories up to Cauchy completion (i.e. up to Morita equivalence), then one may work with $\Map(\QDist)$ in place of the $2$-category $\QCat$. 

\section{The presheaf and copresheaf monads on $\Map(\QDist)$}\label{presheaf and copresheaf monads}

In this section, we extend the constructions of the presheaf and copresheaf monads from $\QCat$ to $\Map(\QDist)$. 

Each left adjoint $\CQ$-distributor $\zeta\colon\bbA\oto\bbB$ gives rise to a pair of $\CQ$-functors
\[\zeta^\ra\colon\CP\bbA\to\CP\bbB,\quad \zeta^\ra(\mu)=\mu\circ\zeta^*, \]
\[\zeta^\la\colon\CP\bbB\to\CP\bbA,\quad\zeta^\la(\lambda)=\lam\circ\zeta.\]

\begin{prop}\label{direct-image}
Let $\zeta\colon\bbA\oto\bbB$ be a left adjoint $\CQ$-distributor. Then $\zeta^\ra\dv\zeta^\la\colon\CP\bbB\to\CP\bbA$ in $\QCat$.
\end{prop}
\begin{proof}
The verification follows by straightforward calculation.
\end{proof}




This yields the \emph{presheaf monad} 
\[(\widehat{\CP},\sm,\iota)\]
on $\Map(\QDist)$ defined as follows:
\begin{itemize}
\item the functor $\widehat{\CP}$ assigns to each $\CQ$-category $\bbA$ the presheaf category $\CP\bbA$, and to each left adjoint $\CQ$-distributor $\zeta\colon\bbA\oto\bbB$ the graph of the $\CQ$-functor $\zeta^\ra\colon\CP\bbA\to\CP\bbB$;
\item for each $\CQ$-category $\bbA$, the component of the unit $\iota_\bbA\colon \bbA\oto \CP\bbA$ is the graph of the Yoneda embedding $\sY_\bbA\colon\bbA\to\CP\bbA$;
\item for each $\CQ$-category $\bbA$, the component of the multiplication $\sm_\bbA\colon\CP\CP\bbA\oto\CP\bbA$ is the graph of the $\CQ$-functor $\sup\nolimits_{\CP\bbA}\colon\CP\CP\bbA\to\CP\bbA$.
\end{itemize}

\begin{prop}\label{T-dual-KZ}
The monad $(\widehat{\CP},\sm,\iota)$ is a co-KZ-monad on the category $\Map(\QDist)$.
\end{prop}
\begin{proof}
Let $\bbA$ be a $\CQ$-category. Since the presheaf category $\CP \bbA$ is cocomplete, 
the $\CQ$-functor  $\sup\nolimits_{\CP \bbA}$ is left adjoint to the Yoneda embedding $\sY_{\CP \bbA}$, that is,
\[\sup\nolimits_{\CP\bbA}\dv\sY_{\CP\bbA}.\]
By Proposition \ref{ad-f-d}\ref{ad-f-d-g}, this adjunction corresponds, in $\QDist$, to
\[\iota_{\widehat{\CP} \bbA}=\iota_{\CP\bbA}=(\sY_{\CP \bbA})_\natural\dashv(\sup\nolimits_{\CP \bbA})_\natural=\sm_\bbA.\]
Hence $(\widehat{\CP},\sm,\iota)$ is a co-KZ-monad.
\end{proof}

With the co-KZ-monad $(\widehat{\CP},\sm,\iota)$ in place, we can now describe its algebras explicitly.  A $\widehat{\CP}$-algebra is a pair $(\bbA,\phi)$ consisting of a $\CQ$-category $\bbA$ together with a left adjoint $\CQ$-distributor $\phi\colon\CP\bbA\oto\bbA$ that is right adjoint to $(\sY_\bbA)_\natural$. In particular, since 
\[(\sY_\bbA)_\natural\dv(\sY_\bbA)^\natural\colon\CP\bbA\oto\bbA,\]and the right adjoint of a morphism in the quantaloid $\QDist$ is unique whenever it exists, it follows that $\phi$ is precisely the cograph of the Yoneda embedding $\sY_\bbA$. 

Before discussing homomorphisms between $\widehat{\CP}$-algebras, we first record a lemma that will be used in the subsequent analysis.

\begin{lem}\label{T-h}
Let $\zeta\colon\bbA\oto\bbB$ be a left adjoint $\CQ$-distributor. 
\begin{enumerate}[label=(\arabic*)]
\item\label{ra-cd} $\zeta$ is a right adjoint $\CQ$-distributor if and only if
\begin{equation}\label{ra-eq}
\Big(\zeta\circ(\sY_\bbA)^\natural\Big)(\mu,y)=\CP\bbA(\mu,\zeta(-,y))
\end{equation}
for any $\mu\in(\CP\bbA)_0$ and $y\in\bbB_0$.
\item\label{la-cd} $\zeta^*$ is a left adjoint $\CQ$-distributor if and only if 
\begin{equation}\label{la-eq}
\Big((\sYd_\bbA)_\natural\circ\zeta^*\Big)(y,\lambda)=\CPd\bbA(\zeta^*(y,-),\lambda)
\end{equation}
for any $\lambda\in(\CPd\bbA)_0$ and $y\in\bbB_0$.
\end{enumerate}
\end{lem}
\begin{proof}
We prove \ref{ra-cd} for example.
Let $\zeta\colon\bbA\oto\bbB$ be a $\CQ$-distributor that is simultaneously a left and a right adjoint. 
For any $\mu\in(\CP\bbA)_0$ and $y\in\bbB_0$, we have
\begin{align*}
\Big(\zeta\circ(\sY_\bbA)^\natural \Big)(\mu,y)&=\zeta(-,y)\circ\CP\bbA(\mu,\sY_\bbA(-))\\
&=\zeta(-,y)\circ(\bbA\swarrow\mu)&(\text{Equation}\ (\ref{sup}))\\
&=\Big(\zeta\circ(\bbA\lda\mu)\Big)(y)\\
&=\Big((\zeta\circ\bbA)\lda\mu\Big)(y) &(\text{Proposition}\ (\ref{Q-dv-cd-3}))\\
&=(\zeta\lda\mu)(y)\\
&=\zeta(-,y)\lda\mu\\
&=\CP\bbA(\mu,\zeta(-,y)).
\end{align*}

Conversely, suppose that equation (\ref{ra-eq}) holds for all $\mu\in(\CP\bbA)_0$ and $y\in \bbB_0$. Then for any $b,b'\in\bbB_0$,
\begin{align*}
\Big(\zeta\circ(\bbA\lda\zeta)\Big)(b,b')&=\zeta(-,b')\circ(\bbA\lda\zeta)(b,-)\\
&=\zeta(-,b')\circ\Big(\bbA(-,-)\lda\zeta(-,b)\Big)\\
&=\zeta(-,b')\circ\CP\bbA(\zeta(-,b),\sY_\bbA(-))\\
&=\zeta(-,b')\circ(\sY_\bbA)^\natural(\zeta(-,b),-)\\
&=\CP\bbA(\zeta(-,b),\zeta(-,b'))&(\text{Equation}\ (\ref{ra-eq}))\\
&=\zeta(-,b')\lda\zeta(-,b)\\
&\geq \bbB(b,b')
\end{align*}
for any $b,b'\in\bbB_0$. Thus $\zeta\colon\bbA\oto\bbB$ is a right adjoint $\CQ$-distributor.
\end{proof}

\begin{prop}\label{Pre-md-al}
Let $\bbA$ and $\bbB$ be $\widehat{\CP}$-algebras. A left adjoint $\CQ$-distributor $\zeta\colon\bbA\oto\bbB$ is the $\widehat{\CP}$-algebra homomorphism if and only if it is  a right adjoint $\CQ$-distributor.
\end{prop}
\begin{proof}
Let $(\bbA,(\sY_\bbA)^\natural)$ and $(\bbB,(\sY_\bbB)^\natural)$ be $\widehat{\CP}$-algebras, and let $\zeta\colon\bbA\oto\bbB$ be a
left adjoint $\CQ$-distributor. Then  
\begin{align*}
\zeta\ \text{is a}\ \widehat{\CP}\text{-algebra homomorphism}&\iff (\sY_\bbB)^\natural\circ(\widehat{\CP}\zeta)=\zeta\circ(\sY_\bbA)^\natural\\
&\iff (\sY_\bbB)^\natural\circ(\zeta^\ra)_\natural=\zeta\circ(\sY_\bbA)^\natural\\
&\iff \CP\bbB(\zeta^\ra(-),\sY_\bbB(-))=\zeta\circ(\sY_\bbA)^\natural.
\end{align*}
Suppose now that $\zeta\colon\bbA\oto\bbB$ is a $\widehat{\CP}$-algebra homomorphism. Then for any $\mu\in(\CP\bbA)_0$ and $y\in\bbB_0$, we have
\begin{align*}
(\zeta\circ(\sY_\bbA)^\natural)(\mu,y)
&=\CP\bbB(\zeta^\ra(\mu),\sY_\bbB(y))\\
&=\CP\bbA(\mu,\zeta^\la(\sY_\bbB(y)))& (\zeta^\ra\dv\zeta^\la)\\
&=\CP\bbA(\mu,\sY_\bbB(y)\circ\zeta)\\
&=\CP\bbA(\mu,\zeta(-,y)).
\end{align*}
Hence, by Lemma \ref{T-h}(\ref{ra-eq}), $\zeta$ is also a right adjoint $\CQ$-distributor. 
The converse implication follows by a similar argument.
\end{proof}

Having constructed the presheaf monad on $\Map(\QDist)$, we now proceed to the dual situation. Each left adjoint $\CQ$-distributor $\zeta\colon\bbA\oto\bbB$ gives rise to a pair of $\CQ$-functors
\[\zeta^{\nla}\colon\CPd\bbB\to\CPd\bbA,\quad \zeta^{\nla}(\lam)=\zeta^*\circ\lam, \]
\[\zeta^{\nra}\colon\CPd\bbA\to\CPd\bbB,\quad \zeta^{\nra}(\mu)=\zeta\circ\mu.\]
By the same reasoning as in Proposition \ref{direct-image}, one has an adjunction $\zeta^{\nla}\dv\zeta^{\nra}\colon\CPd\bbA\to\CPd\bbB$. This yields the \emph{copresheaf monad} 
\[(\widehat{\CPd},\sm^\dag,\iota^\dag)\]
on $\Map(\QDist)$ defined as follows:
\begin{itemize}
\item the functor $\widehat{\CPd}$ sends each $\CQ$-category $\bbA$ to $\CPd\bbA$, and sends each left adjoint $\CQ$-distributor $\zeta\colon\bbA\oto\bbB$ to $(\zeta^{\nra})_\natural\colon\CPd\bbA\oto\CPd\bbB$;
\item for each $\CQ$-category $\bbA$, the component of the unit $\iota^\dag_\bbA\colon\bbA\oto\CPd\bbA$ is the graph of the co-Yoneda embedding $\sYd_\bbA\colon\bbA\to\CPd\bbA$;
\item for each $\CQ$-category $\bbA$, the component of the multiplication $\sm^\dag_\bbA\colon\CPd\CPd\bbA\oto\CPd\bbA$ is the graph of the $\CQ$-functor $\inf\nolimits_{\CPd\bbA}\colon\CPd\CPd\bbA\to\CPd\bbA$.
\end{itemize}
Dually, the copresheaf monad $(\widehat{\CPd},\sm^\dag,\iota^\dag)$ is a KZ-monad on $\Map(\QDist)$. A $\widehat{\CPd}$-algebra is pair $(\bbA,\phi)$ consisting of a $\CQ$-category $\bbA$ and a $\CQ$-distributor $\phi\colon\CPd\bbA\to\bbA$ that is left adjoint to $(\sYd_\bbA)_\natural$.

\begin{lem}\label{la-cd-2}
Let $(\bbA,\phi)$ and $(\bbB,\psi)$ be $\widehat{\CPd}$-algebras and $\zeta\colon\bbA\oto\bbB$ a left adjoint $\CQ$-distributor. $\zeta^*$ is a left adjoint $\CQ$-distributor if and only if 
\begin{equation}\label{la-cd-eq-2}
\zeta\circ\phi=\bbB\swarrow\big(\CPd\bbA(\zeta^*(-,-),-)\big).
\end{equation}
\end{lem}
\begin{proof}
For sufficiency, suppose that the equation (\ref{la-cd-eq-2}) holds. Since $(\bbA,\phi)$ is the  $\widehat{\CPd}$-algebra, we have $\phi\dv(\sYd_\bbA)_\natural$ and moreover  $(\sYd_\bbA)_\natural\circ\zeta^*$ is a right adjoint $\CQ$-distributor, that is, 
\[\zeta\circ\phi\dv(\sYd_\bbA)_\natural\circ\zeta^*.\]
It then follows that
\begin{align*}
& \big((\sYd_\bbA)_\natural\circ\zeta^* \big)\swarrow\CPd\bbA(\zeta^*(-,-),-)\\
\geq {}&\Big(\big((\sYd_\bbA)_\natural\circ\zeta^* \big)\swarrow\bbB\Big)\circ\Big(\bbB\swarrow\big(\CPd\bbA(\zeta^*(-,-),-)\big) \Big)\\
={}&\big((\sYd_\bbA)_\natural\circ\zeta^* \big)\circ\big(\zeta\circ\phi\big)&(\text{Equation}\ (\ref{la-cd-eq-2})) \\
\geq{}&\CPd\bbA.&(\zeta\circ\phi\dv(\sY_\bbA)_\natural\circ\zeta^*)
\end{align*}
Consequently, we have
\[(\sYd_\bbA)_\natural\circ\zeta^*\geq\CPd\bbA(\zeta^*(-,-),-).\]
Note that for any $y\in\bbB_0$ and $\mu\in(\CPd\bbA)_0$,
\begin{align*}
\Big((\sYd_\bbA)_\natural\circ\zeta^*\Big)(y,\mu)&=\CPd\bbA(\sYd_\bbA(-),\mu)\circ\zeta^*(y,-)\\
&=\CPd\bbA(\sYd_\bbA(-),\mu)\circ\CPd\bbA(\zeta^*(-,-),\sYd_\bbA(-))\\
&=\Big((\sYd_\bbA)_\natural\circ(\sYd_\bbA)^\natural\Big)(\zeta^*(y,-),\mu). 
\end{align*}
Hence, the equation (\ref{la-eq}) holds, since the inequality 
\[(\sYd_\bbA)_\natural\circ\zeta^*=\big((\sYd_\bbA)_\natural\circ(\sYd_\bbA)^\natural\big)(\zeta^*(-,-),-)\leq\CPd\bbA(\zeta^*(-,-),-)\]
is immediate. By Lemma \ref{T-h}\ref{la-cd}, $\zeta^*$ is therefore a left adjoint $\CQ$-distributor.

For necessity, assume that $\zeta^*$ is a left adjoint $\CQ$-distributor. By Lemma \ref{T-h}\ref{la-cd}, one has 
\[\bbB\lda\big((\sY_\bbA)_\natural\circ\zeta^*\big)=\bbB\lda\big(\CPd\bbA(\zeta^*(-,-),-)\big).\]
Since $\zeta\circ\phi\dv(\sYd_\bbA)_\natural\circ\zeta^*$, it follows that
\[\zeta\circ\phi=\bbB\lda\big((\sY_\bbA)_\natural\circ\zeta^*\big),\]
and hence
\[\zeta\circ\phi=\bbB\lda\big(\CPd\bbA(\zeta^*(-,-),-)\big).\]
\end{proof}

\begin{prop}\label{coPre-md-al}
Let $\bbA$ and $\bbB$ be $\widehat{\CPd}$-algebras. A left adjoint $\CQ$-distributor $\zeta\colon\bbA\oto\bbB$ is the $\widehat{\CPd}$-algebra homomorphism if and only if $\zeta^*$ is a left adjoint.
\end{prop}
\begin{proof}
Let $(\bbA,\phi)$ and $(\bbB,\psi)$ be the $\widehat{\CPd}$-algebras, and let $\zeta\colon\bbA\oto\bbB$ be a left adjoint $\CQ$-distributor. Then  
\begin{align*}
\zeta\ \text{is a}\  \widehat{\CPd}\text{-algebra homomorphism} &\iff \psi\circ(\widehat{\CPd}\zeta)=\zeta\circ\phi\\
&\iff \psi\circ(\zeta^{\nra})_\natural=\zeta\circ\phi\\
&\iff \psi(\zeta^{\nra}(-),-)=\zeta\circ\phi.
\end{align*}
It follows that
\[\psi(\zeta^{\nra}(-),-)=\bbB\lda\big(\CPd\bbA(\zeta^*(-,-),-)\big),\]
since 
\begin{align*}
\psi(\zeta^{\nra}(\mu),y)&=\Big(\bbB\swarrow(\sYd_\bbB)_\natural\Big)(\zeta^{\nra}(\mu),y)&(\psi\dv(\sYd_\bbB)_\natural)\\
&=\bbB(-,y)\swarrow(\sYd_\bbB)_\natural(-,\zeta^{\nra}(\mu))\\
&=\bbB(-,y)\swarrow\CPd\bbB(\sYd_\bbB(-),\zeta^{\nra}(\mu))\\
&=\bbB(-,y)\swarrow\CPd\bbA(\zeta^{\nla}(\sYd_\bbB(-)),\mu)&(\zeta^{\nla}\dv\zeta^{\nra})\\
&=\bbB(-,y)\swarrow\CPd\bbA(\zeta^*\circ\sYd_\bbB(-),\mu)\\
&=\bbB(-,y)\swarrow\CPd\bbA(\zeta^*(-,-),\mu)
\end{align*}
for any $\mu\in(\CPd\bbA)_0$, $y\in\bbB_0$.
Thus $\zeta$ is a $\widehat{\CPd}$-algebra homomorphism if and only if $\bbB\lda\big(\CPd\bbA(\zeta^*(-,-),-)\big)=\zeta\circ\phi $.
By Lemma \ref{la-cd-2}, it follows that $\zeta$ is a $\widehat{\CPd}$-algebra homomorphism if and only if $\zeta^*$ is a left adjoint $\CQ$-distributor.
\end{proof}

\section{Completeness of $\CQ$-categories up to Morita equivalence}\label{M-completeness}

After constructing the presheaf and the copresheaf monads on $\Map(\QDist)$, we are ready to introduce the notions of completeness and cocompleteness up to Morita equivalence.

\begin{defn}
Let $\bbA$ be a $\CQ$-category.
\begin{enumerate}[label=(\arabic*)]
\item $\bbA$ is \emph{cocomplete up to Morita equivalence}, or \emph{$\CM$-cocomplete} for short, if the cograph $(\sY_\bbA)^\natural\colon\CP\bbA\oto\bbA$ of the Yoneda embedding is a left adjoint $\CQ$-distributor.
\item $\bbA$ is \emph{complete up to Morita equivalence}, or \emph{$\CM$-complete} for short, if the graph $(\sYd_\bbA)_\natural\colon\bbA\oto\CPd\bbA$ of the co-Yoneda embedding is a right adjoint $\CQ$-distributor.
\end{enumerate}
\end{defn}

We recall that a $\CQ$-category is cocomplete if and only if it is complete \cite{Stubbe2005}. Under our definitions, one obtains an analogous equivalence between cocompleteness and completeness up to Morita equivalence. Before proving this result, we observe that for any $\CQ$-category $\bbA$, there is an adjunction
\[\ub\dv\lb\colon\CPd\bbA\to\CP\bbA\]
in $\QCat$, defined by
\[\ub(\mu)=\bbA\lda\mu\quad\text{and}\quad\lb(\lam)=\lam\rda\bbA. \]
For every object $a\in\bbA_0$ and for each $\mu\in(\CP\bbA)_0$, $\lam\in(\CPd\bbA)_0$, their evaluations are given by
\begin{equation}\label{ub-Y}
(\ub(\mu))(a)=\CP\bbA(\mu,\bbA(-,a))=(\sY_\bbA)^\natural(\mu,a)   
\end{equation}
and
\begin{equation}\label{lb-Yd}
(\lb(\lam))(a)=\CPd\bbA(\bbA(a,-),\lam)=(\sYd_\bbA)_\natural(a,\lam). 
\end{equation}

\begin{thm}\label{QMC-QMCC-eq}
A $\CQ$-category is $\CM$-cocomplete if and only if it is $\CM$-complete.
\end{thm}
\begin{proof}
Let $\bbA$ be a $\CM$-cocomplete $\CQ$-category. Define a $\CQ$-distributor $\psi$ by
\[\psi:=(\sY_\bbA)^\natural\circ(\lb)_\natural\colon \CPd\bbA\oto\CP\bbA\oto\bbA,\] 
that is, 
\[\psi=\CP\bbA(-\rda\bbA,\sY_\bbA(-))=\bbA\lda(-\rda\bbA).\]
Since $\bbA$ is $\CM$-cocomplete, $(\sY_\bbA)^\natural$ is a left adjoint $\CQ$-distributor, and therefore so is $\psi$. To establish the adjunction $\psi\dv(\sYd_\bbA)_\natural$, it suffices to verify that
\[(\sYd_\bbA)_\natural=\psi^*=\psi\rda\bbA.\]
Indeed, for each $a\in\bbA_0$ and $\lam\in(\CPd\bbA)_0$, we have
\begin{align*}
(\psi\rda\bbA)(a,\lambda)&=\bw\limits_{x\in{\bbA_0}}\psi(\lambda,x)\rda\bbA(a,x)\\
&=\bw\limits_{x\in\bbA_0}\big(\bbA\lda(\lam\rda\bbA)\big)(x)\rda\bbA(a,x)\\
&=(\lam\rda\bbA)(a)&(\ub\dv\lb)\\
&=\CPd\bbA(\bbA(a,-),\lam)\\
&=(\sYd_\bbA)_\natural(a,\lam)&(\text{Equation}\ (\ref{lb-Yd}))
\end{align*}
Consequently, $\psi\dv(\sYd_\bbA)_\natural$, and thus $\bbA$ is $\CM$-complete. The converse implication follows by a similar reasoning.
\end{proof}

\begin{rem}
For every $\CQ$-category $\bbA$, combining the isomorphism (\ref{Q-Dist-op}) with $\CPd\bbA\cong(\CP(\bbA^\op))^\op$, one obtains that $\bbA$ is $\CM$-complete if and only if $\bbA^\op$ is $\CM$-cocomplete.
\end{rem}

\begin{defn}
Let $\bbA$ and $\bbB$ be $\CM$-cocomplete $\CQ$-categories and $\zeta\colon\bbA\to\bbB$ be a left adjoint $\CQ$-distributor.
\begin{enumerate}[label=(\arabic*)]
\item $\zeta$ is \emph{$\CM$-cocontinuous} if it is also a right adjoint $\CQ$-distributor. 
\item $\zeta$ is \emph{$\CM$-continuous} if $\zeta^*$ is a left adjoint $\CQ$-distributor.
\end{enumerate}
\end{defn}

It is well known that $\QSup$ and $\QInf$ are the categories of algebras for the presheaf and copresheaf monads on $\QCat$, respectively. Accordingly, we denote by
\[\QMCC:=\Map(\QDist)^{\widehat{\CP}}\]
the category of $\CM$-cocomplete $\CQ$-categories and $\CM$-cocontinuous left adjoint $\CQ$-distributors. Similarly, we denote by
\[\QMC:=\Map(\QDist)^{\widehat{\CPd}}\]
the category of $\CM$-complete $\CQ$-categories and $\CM$-continuous left adjoint $\CQ$-distributors. 

\begin{thm}\label{Q-Sup-Map(Q-Dist)}
Let $\zeta\colon\bbA\oto\bbB$ be a left adjoint $\CQ$-distributor with $\bbB$ $\CM$-cocomplete. There exists a unique $\CM$-cocontinuous left adjoint $\CQ$-distributor $\eta\colon\CP\bbA\oto\bbB$ such that
\[\zeta=\eta\circ(\sY_\bbA)_\natural.\]
\end{thm}
\begin{proof}
\textbf{Existence.} Since $\bbB$ is $\CM$-cocomplete, the cograph of the Yoneda embedding $\sY_\bbB$ is a left adjoint $\CQ$-distributor. We define $\eta$ as the composite of left adjoint $\CQ$-distributors:
\[(\sY_\bbB)^\natural\circ(\zeta^\ra)_\natural\colon\CP\bbA\oto\CP\bbB\oto\bbB,\]
that is,
\[\eta=\CP\bbB(\zeta^\ra(-),\sY_\bbB(-)).\]
For any $x\in\bbA_0$ and  $y\in\bbB_0$, we have
\begin{align*}
\Big(\eta\circ(\sY_\bbA)_\natural \Big)(x,y)&=\eta(\sY_\bbA(x),y)\\
&=\CP\bbB(\zeta^\ra(\sY_\bbA(x)),\sY_\bbB(y))\\
&=\CP\bbA(\sY_\bbA(x),\zeta^\la(\sY_\bbB(y)))&(\zeta^\ra\dv\zeta^\la)\\
&=\CP\bbA(\sY_\bbA(x),\zeta(-,y))\\
&=\zeta(x,y).
\end{align*}
Thus, $\eta\circ(\sY_\bbA)_\natural=\zeta$. It remains to show that $\eta$ is $\CM$-cocontinuous. By Proposition \ref{ad-f-d} and \ref{direct-image}, we have $(\zeta^\la)_\natural\dv(\zeta^\ra)_\natural$. Hence, $\eta$ is also a right adjoint $\CQ$-distributor.

\textbf{Uniqueness.} Suppose that $\eta'\colon\CP\bbA\oto\bbB$ is another $\CM$-cocontinuous left adjoint $\CQ$-distributor such that $\zeta=\eta'\circ(\sY_\bbA)_\natural$. For any $\Phi\in(\CP\CP\bbA)_0$ and $y\in\bbB_0$, we have
\[\eta'(\sup\nolimits_{\CP\bbA}\Phi,y)=\CP\CP\bbA(\Phi,\eta'(-,y))\]
since
\begin{align*}
\Big(\eta'\circ(\sY_{\CP\bbA})^\natural \Big)(\Phi,y)&=\eta'(-,y)\circ\CP\CP\bbA(\Phi,\sY_{\CP\bbA}(-))\\
&=\eta'(-,y)\circ\CP\bbA(\sup\nolimits_{\CP\bbA}\Phi,-)&(\sup\nolimits_{\CP\bbA}\dv\sY_{\CP\bbA})\\
&=\eta'(\sup\nolimits_{\CP\bbA}\Phi,y)\\
&=\CP\CP\bbA(\Phi,\eta'(-,y))&(\eta'\text{ is }\CM\text{-cocontinuous})
\end{align*}
Observe that every presheaf $\mu$ on $\bbA$ is the supremum of $\sY_\bbA^\ra(\mu)$. Then,
\begin{align*}
\eta'(\mu,y)&=\eta'(\sup\nolimits_{\CP\bbA}(\sY_\bbA^\ra)(\mu),y)\\
&=\CP\CP\bbA(\sY_\bbA^\ra(\mu),\eta'(-,y))\\
&=\CP\bbA(\mu,\sY_\bbA^\la(\eta'(-,y)))&(\sY_\bbA^\ra\dv\sY_\bbA^\la)\\
&=\CP\bbA(\mu,(\eta'\circ(\sY_\bbA)_\natural)(-,y))\\
&=\CP\bbA(\mu,\zeta(-,y))&(\eta'\circ(\sY_\bbA)_\natural=\zeta)\\
&=\CP\bbA(\mu,\sY_\bbB(y)\circ\zeta)\\
&=\CP\bbA(\mu,\zeta^\la(\sY_\bbB(y)))\\
&=\CP\bbB(\zeta^\ra(\mu),\sY_\bbB(y))&(\zeta^\ra\dv\zeta^\la)\\
&=\Big((\sY_\bbB)^\natural\circ(\zeta^\ra)_\natural\Big)(\mu,y)
\end{align*}
for all $\mu\in(\CP\bbA)_0$ and $y\in\bbB_0$.
Therefore,
\[\eta'=(\sY_\bbB)^\natural\circ(\zeta^\ra)_\natural,\]
and thus $\eta'=\eta$, as required.

\end{proof}

\begin{rem}
The above theorem, in fact, establishes an adjunction:
\begin{equation} 
\bfig
\morphism|a|/@{->}@<6pt>/<1000,0>[\Map(\QDist)`\QMCC;\widehat{\CP}]
\morphism(1000,0)|b|/@{->}@<6pt>/<-1000,0>[\QMCC`\Map(\QDist);U_{\QMCC}]
\place(550,0)[\bot]
\efig
\end{equation}
where $U_{\QMCC}$ denotes the forgetful functor.
The monad induced by this adjunction is precisely the presheaf on $\Map(\QDist)$.

In particular, for any $\CM$-cocontinuous left adjoint $\CQ$-distributor $\zeta\colon\bbA\oto\bbB$ with $\bbB$ cocomplete, there exists a unique left adjoint $\CQ$-functor $F\colon\CP\bbA\to\bbB$ such that
\[\zeta=(F\circ\sY_\bbA)_\natural.\]
Restricting  $(-)_\natural\colon\QCat\to\Map(\QDist)$ to the categories $\QSup$ and $\QMCC$, it follows that every $\CM$-cocomplete $\CQ$-category $\bbA$ admits a reflection along \[(-)_\natural\colon\QSup\to\QMCC,\]
given by $(\CP\bbA,(\sY_\bbA)_\natural)$.
\end{rem}


We now present several examples of $\CM$-cocomplete $\CQ$-categories which are not $\CQ$-complete, showing that $\CM$-cocompleteness is strictly weaker than $\CQ$-completeness.

\begin{exmp}
\begin{enumerate}[label=(\arabic*)]
\item The frame $F=\{\bot,p,q,k\}$, illustrated by the Hasse diagram
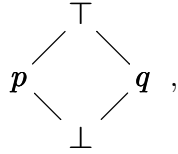
\[\bfig
\Atriangle/-`-`/<200,200>[\top`p`q;``]
\Vtriangle(0,-200)/`-`-/<200,200>[p`q`\bot;``]
\place(500,0)[\text{,}]
\efig\]
is a one-obejct quantaloid.
Let $X=\{x,y\}$ be an $F$-category with 
\[X(x,y)=p\quad\text{and}\quad X(y,x)=q.\]
The corresponding values of $(\sY_X)^\natural\colon\CP X\to X$ and $((\sY_X)^\natural\rda X)^\op\colon X^\op\to(\CP X)^\op$ are listed in the table:
\[\begin{array}{c|cc}
	    (\sY_X)^\natural&x&y\\
	\hline \mu(x)=\mu(y)=\bot&k&k\\
    \mu(x)=\bot,\ \mu(y)=q&k&k\\
    \mu(x)=p,\ \mu(y)=\bot& k&k\\
    \mu(x)=\mu(y)=p&q&k\\
    \mu(x)=p,\ \mu(y)=q&k&k\\
    \mu(x)=p,\ \mu(y)=k&q&k\\
    \mu(x)=\mu(y)=q&k&p\\
    \mu(x)=k,\ \mu(y)=q&k&p\\
    \mu(x)=\mu(y)=k&q&p\
\end{array}\quad\text{and}\quad\begin{array}{c|cc}
	    ((\sY_X)^\natural\rda X)^\op&x&y\\
	\hline \mu(x)=\mu(y)=\bot&p&q\\
    \mu(x)=\bot,\ \mu(y)=q&p&q\\
    \mu(x)=p,\ \mu(y)=\bot& p&q\\
    \mu(x)=\mu(y)=p&q&k\\
    \mu(x)=p,\ \mu(y)=q&p&q\\
    \mu(x)=p,\ \mu(y)=k&p&k\\
    \mu(x)=\mu(y)=q&k&q\\
    \mu(x)=k,\ \mu(y)=q&k&q\\
    \mu(x)=\mu(y)=k&k&k\
\end{array} \]
From simple calculations it follows that 
\[\sY_X^\natural\colon\CP X\oto X\]
is a left adjoint $F$-distributor. Hence $X$ is $\CM$-cocomplete. Moreover, $X$ is not cocomplete.
\item Let $[-\infty,\infty]$ be the extended real line equipped with the order ``$\geq$''. Then $([-\infty,\infty],+,0)$ is a one-object quantaloid. Let $X=\{x,y\}$ be a $[-\infty,+\infty]$-category with $X(x,y)=-1$, $X(y,x)=1$ and $X(x,x)=X(y,y)=0$. 
Given any presheaf $\mu\colon X\oto *$, we obtain $\mu(y)=\mu(x)+1$ from
\[\mu(x)+X(y,x)\geq \mu(y)\quad \text{and}\quad \mu(y)+X(x,y)\geq \mu(x).\]
Then we have
\begin{align*}
(\sY_X)^\natural(\mu,x)&=\bv\limits_{a\in X}\mu(a)\ra X(a,x)\\
&=\Big(X(x,x)-\mu(x) \Big)\vee\Big(X(y,x)-\mu(y) \Big)\\
&=\Big(0-\mu(x)\Big)\vee\Big(1-(1+\mu(x))\Big)\\
&=-\mu(x)
\end{align*}
and
\begin{align*}
(\sY_X)^\natural(\mu,y)&=\bv\limits_{a\in X}\mu(a)\ra X(a,y)\\
&=\Big(X(x,y)-\mu(x) \Big)\vee\Big(X(y,y)-\mu(y) \Big)\\
&=\Big(-1-\mu(x)\Big)\vee\Big(0-(1+\mu(x))\Big)\\
&=-1-\mu(x).
\end{align*}
Hence
\[((\sY_X)^\natural\rda X)(x,\mu)=\mu(x)\quad\text{and}\quad((\sY_X)^\natural\rda X)(y,\mu)=1+\mu(x).\]
It follows immediately that
\[\CP X(\mu,\mu)\leq \big(((\sY_X)^\natural\rda X)\circ(\sY_X)^\natural\big)(\mu,\mu).\] 
Since $\mu$ is arbitrary, $\sY_X^\natural\colon\CP X\oto X$ is a left adjoint $[-\infty,+\infty]$-distributor. Consequently, $X$ is $\CM$-cocomplete, although it is not cocomplete.
\end{enumerate}
\end{exmp}

Recall that a $\CQ$-category is complete if and only if it is tensored and conically cocomplete. In the following, we examine the corresponding notions in the context of Morita equivalence.

\begin{defn}
Let $\bbA$ be a $\CQ$-category.
\begin{enumerate}[label=(\arabic*)]
\item $\bbA$ is \emph{$\CM$-tensored} if for every $a\in\bbA_0$,  the cograph 
\[(\bbA(a,-))^\natural\colon \CP|a|\oto\bbA\]
is a left adjoint $\CQ$-distributor. 
\item $\bbA$ is \emph{$\CM$-cotensored} if for every $a\in\bbA_0$, the graph 
\[(\bbA(-,a))_\natural\colon\bbA \oto\CPd|a|\] is a right adjoint $\CQ$-distributor.
\end{enumerate}
\end{defn}

\begin{prop} \label{cp-tensor}
Let $\bbA$ be a $\CQ$-category.
\begin{enumerate}[label=(\arabic*)]
\item\label{cp-tensor-1} For all $a\in\bbA_0$ and $f\in{\CP|a|}_0$,
\[(\bbA(a,-))^\natural(f,-)=(\sY_\bbA)^\natural(f\circ\sY_\bbA(a),-).\]
\item\label{cp-tensor-2} For all $a\in\bbA_0$ and $g\in{\CPd|a|}_0$,
\[(\bbA(-,a))_\natural(-,g)=(\sY_\bbA)_\natural(-,\sYd_\bbA(a)\circ g).\]
\end{enumerate}
\end{prop}
\begin{proof}
We verify \ref{cp-tensor-1}; the proof of \ref{cp-tensor-2} follows by the similar argument and will be omitted. Let $a\in\bbA_0$ and $f\in{\CP|a|}_0$. Then we have
\begin{align*}
(\bbA(a,-))^\natural(f,-)&=\CP(\left| a\right|)(f,\bbA(a,-))\\
&=\bbA(a,-)\lda f\\
&=(\bbA\lda\bbA(-,a))\lda f\\
&=\bbA\lda(f\circ\sY_\bbA(a))\\
&=\CP\bbA(f\circ\sY_\bbA(a),\sY_\bbA(-))\\
&=(\sY_\bbA)^\natural(f\circ\sY_\bbA(a),-).
\end{align*}
\end{proof}

\begin{defn}\label{M-conical-df}
Let $\bbA$ be a $\CQ$-category.
\begin{enumerate}[label=(\arabic*)]
\item\label{M-conical-df-1} $\bbA$ is \emph{$\CM$-conically cocomplete} if, for any subsets $\{x_i\}\subseteq\bbA_0$ and $\{f_i\colon *_{|x_i|}\oto *_X\}$ with $X\in \CQ_0$, the $\CQ$-distributor
\[\bw\limits_{i}\Big(\big(\bbA(x_i,-) \big)^\natural(f_i,-)\Big)\colon *_X\oto\bbA\]
is a left adjoint $\CQ$-distributor whenever each individual $\big(\bbA(x_i,-) \big)^\natural(f_i,-)$ is a left adjoint.
\item \label{M-conical-df-2}$\bbA$ is \emph{$\CM$-conically complete} if, for any subsets $\{x_i\}\subseteq\bbA$ and $\{g_i\colon*_X\oto *_{|x_i|}\}$ with $X\in\CQ_0$, the $\CQ$-distributor
\[\bw\limits_i\Big(\big(\bbA(-,x_i)\big)_\natural\Big)(-,g_i)\colon\bbA\oto *_X\] 
is a right adjoint whenever each individual $\big(\bbA(-,x_i)\big)_\natural(-,g_i)$ is a right adjoint.
\end{enumerate}
\end{defn}

\begin{rem}
When the $\CQ$-category $\bbA$ is Cauchy complete, it is straightforward to verify that:
\begin{itemize}
\item $\bbA$ is $\CM$-(co)complete if and only if it is (co)complete;
\item $\bbA$ is $\CM$-(co)tensored if and only if it is (co)tensored.
\end{itemize} 

However, it is less immediate that $\bbA$ is $\CM$-conically (co)complete if and only if it is conically (co)complete. Let $\bbA$ be a $\CM$-conically cocomplete $\CQ$-category, and let $\{x_i\}_{i\in I}\subseteq\bbA_X$.
For each $i\in I$ and $a\in\bbA_0$, 
\begin{align*}
\big(\bbA(x_i,-)\big)^\natural(1_X,a)&=\CP X(1_X,\bbA(x_i,a))\\
&=\bbA(x_i,a)\swarrow 1_X\\
&=\bbA(x_i,a),
\end{align*}
so each $\bbA(x_i,-)=\big(\bbA(x_i,-)\big)^\natural(1_X,-)$ is a left adjoint $\CQ$-distributor. By Definition \ref{M-conical-df}, 
\[\bw\limits_{i}\Big(\big(\bbA(x_i,-)\big)^\natural(1_X,-)\Big)\]
is therefore also a left adjoint $\CQ$-distributor.
Moreover, we have
\[\bw\limits_{i}\Big(\big(\bbA(x_i,-)\big)^\natural(1_X,-)\Big)=(\sY_\bbA)^\natural(\bv\limits_i \sY_\bbA(x_i),-)\]
since
\begin{align*}
\bw\limits_{i}\Big(\big(\bbA(x_i,-)\big)^\natural(1_X,a)\Big)&=\bw\limits_i \CP X(1_X,\bbA(x_i,a))\\
&=\bw\limits_i \bbA(x_i,-)\\
&=\bw\limits_i \Big(\bbA(-,a)\swarrow\bbA(-,x_i)\Big)\\
&=\bbA(-,a)\swarrow\big(\bv\limits_i \bbA(-,x_i)\big)\\
&=\CP\bbA(\bv\limits_i\bbA(-,x_i),\bbA(-,a))\\
&=\CP\bbA(\bv\limits_i\sY_\bbA(x_i),\sY_\bbA(a))\\
&=(\sY_\bbA)^\natural(\bv\limits_i\sY_\bbA(x_i),a)
\end{align*}
for all $a\in\bbA_0$.
Hence, the $\CQ$-distributor 
\[(\sY_\bbA)^\natural(\bv\limits_i\sY_\bbA(x_i),-)\]
is a left adjoint.
Dually, if $\bbA$ is a $\CM$-conically complete $\CQ$-category and $\{x_i\}_{i\in I}\subseteq\bbA_X$, then the $\CQ$-distributor 
\[(\sYd_\bbA)_\natural(-,\bv\limits_i\sYd_\bbA(x_i))\]
is a right adjoint. 
In this way, $\CM$-conical (co)completeness coincides with conical (co)completeness for $\CQ$-categories studied in \cite{Stubbe2006}.
\end{rem}

In what follows, we examine whether $\CM$-cocompleteness is equivalent to being both $\CM$-tensored and $\CM$-conically cocomplete, as in the classical case. 

\begin{thm}
A $\CQ$-category $\bbA$ is $\CM$-cocomplete if, and only if it is $\CM$-tensored and $\CM$-conically cocomplete.
\end{thm}
\begin{proof}
Let $\bbA$ be a $\CQ$-category. Suppose that $\bbA$ is both $\CM$-tensored and $\CM$-conically cocomplete. Then, for any $\mu\in(\CP\bbA)_0$, we have
\begin{equation}\label{conical}
\bw\limits_{a\in\bbA_0}(\bbA(a,-))^\natural(\mu(a),-)=(\sY_\bbA)^\natural(\mu,-)
\end{equation}
since
\begin{align*}
\bw\limits_{a\in\bbA_0}(\bbA(a,-))^\natural(\mu(a),-)&=\bw\limits_{a\in\bbA_0}\CP(\left| a\right|)(\mu(a),\bbA(a,-))\\
&=(\bbA\lda\mu)(-)\\
&=\CP\bbA(\mu,\sY_\bbA(-))\\
&=(\sY_\bbA)^\natural(\mu,-).
\end{align*}
Therefore, by the definitions of $\CM$-tensored and $\CM$-conically cocomplete $\CQ$-categories, $(\sY_\bbA)^\natural\colon\CP\bbA\oto\bbA$ is a left adjoint $\CQ$-distributor. 

Conversely, suppose that $\bbA$ is $\CM$-cocomplete, i.e.  $(\sY_\bbA)^\natural$ is a left adjoint $\CQ$-distributor. By Proposition \ref{cp-tensor}, for each $a\in\bbA_0$, the $\CQ$-distributor $\big(\bbA(a,-) \big)^\natural\colon\CP|a|\oto\bbA$ is a left adjoint. Thus $\bbA$ is $\CM$-tensored. Next, let $\{x_i\}\subseteq\bbA_0$ and $\{f_i\colon*_{|x_i|}\oto*_X\}$ with $X\in\CQ_0$. From the proof of Proposition \ref{cp-tensor}, we have
\begin{align*}
\bw\limits_i\big(\bbA(x_i,-)\big)^\natural(f_i,-)&=\bw\limits_i\big(\bbA\lda(f_i\circ\sY_\bbA(x_i))\big)\\
&=\bbA\lda\big(\bv\limits_i f_i\circ\sY_\bbA(x_i)\big)\\
&=(\sY_\bbA)^\natural(\bv\limits_i f_i\circ\sY_\bbA(x_i),-).
\end{align*}
Thus the $\CQ$-distributor
\[\bw\limits_i\Big(\big(\bbA(x_i,-)\big)^\natural(f_i,-)\Big)\colon *_X\oto\bbA\]
is also a left adjoint. 
Consequently, $\bbA$ is $\CM$-conically cocomplete.
\end{proof}



\bibliographystyle{abbrv}
\bibliography{lili}

\end{document}